\input xy
\xyoption{all}

\CompileMatrices
\magnification=1200 \hsize=13.3cm \hoffset=0cm

\footline={\vbox{\medskip\hrule\medskip\line
{\pc Chapoton-Livernet\hfill\tenrm\folio}}}


\def\date{\the\day\ \ifcase\month\or janvier \or
f\'evrier\or mars\or avril\or mai\or juin\or juillet\or aout\or septembre
\or octobre\or novembre\or d\'ecembre\fi\ \the\year}
\def\today{\the\month-\the\day-\the\year}

\def\mk{\par\vskip 10mm}



\def\bla{\vskip 21pt plus 5pt minus 5pt}

\def\tibla{\vskip 10pt plus 2pt minus 2pt}

\def\mk{\par\vskip 10mm}
\def\ilestpreferabledecouperici{\penalty -100}


 \newcount\parno
\newcount\secno
\newcount\chapno
\newcount\enonceno
\newcount\eqnrefno
\newcount\figrefno
\parno=0
\secno=0
\chapno=0
\enonceno=0
\eqnrefno=0
\figrefno=0


\def\leaderfill{{\leaders\hbox to 1em{\hss. \hss}\hfill}} 

 
\def\chap#1{
\global\advance\chapno by 1
\addchaptocontents{#1}
\secno=0
\enonceno=0
\eqnrefno=0
\figrefno=0
\mk
\ilestpreferabledecouperici
\centerline{\twelvebf\uppercase\expandafter{\romannumeral\the\chapno}.\ #1}
\nobreak}

\def\initappendice{
\chapno=64
\addinitapptocontents
\global\let\initappendice=\relax}

\def\appendice#1{
\initappendice
\global\advance\chapno by 1
\addchaptocontents{#1}
\secno=0
\enonceno=0
\eqnrefno=0
\figrefno=0
\mk
\ilestpreferabledecouperici
\centerline{\bf\char\chapno\ #1}
\nobreak}

\def\section#1{
\bla
\global\advance\secno by 1
\addsectocontents{#1}
\parno=0
\enonceno=0
\eqnrefno=0
\figrefno=0
\centerline{\twelvebf\the\secno. #1}
\par\nobreak\bla}

\def\paragr#1{
\enonceno=0
\eqnrefno=0
\figrefno=0
\global\advance\parno by 1
\addpartocontents{#1}
\tibla
\ilestpreferabledecouperici
{\bf\the\secno.\the\parno. #1}
\par\nobreak\tibla}

\def\biblio#1{
\vskip 0.7cm
\line{\hfill\twelvebf #1\hfill}
\addbibtocontents{#1}
\par\nobreak
\vskip 0.5cm}

\def\intro#1{
\vskip 1cm
\line{\hfill\twelvebf #1\hfill}
\addintrotocontents{#1}
\par\nobreak
\vskip 1cm}


\def\enonce#1{
\advance \enonceno by 1
\ifnum\parno=0
\def\id{\the\secno .\the\enonceno .}
\else
\def\id{\the\secno .\the\parno .\the\enonceno .}
\fi
\tibla
\ilestpreferabledecouperici
\nd{\pc\id\  #1 }}

\def\debut{\begingroup\it}
\def\fin{\endgroup}


\def\DEF{\enonce{Definition.}}

\def\TH{\enonce{Theorem.}\debut }
\def\PROP{\enonce{Proposition.}\debut }
\def\LEM{\enonce{Lemma.}\debut}

\def\COR{\enonce{Corollary.} \debut}

\def\tit#1{\enonce{#1.}}

\def\dem{
\tibla
\ilestpreferabledecouperici
\nd{\sl Proof.} }

\def\ssnumero#1#2{\expandafter\def\csname#1\endcsname{\nd{\bf #2}. }}
\ssnumero{lemme}{Lemma}
\ssnumero{nota}{Notation}
\ssnumero{contenu}{Contenu du texte}

%
%


\def\ifundefined#1{\expandafter\ifx\csname #1\endcsname\relax }


\newwrite\auxfile
\immediate\openout\auxfile=\jobname.aux

\newtoks\maquereau
\def\ecrire#1#2{\relax
\maquereau=\expandafter{\csname #1\endcsname}\relax
\immediate\write\auxfile{ \def\the\maquereau{#2}}}


\def\label#1{\relax 
\ecrire{#1chapno}{\the\chapno}\relax 
\ifnum\enonceno=0\ifnum\parno=0\relax
\ecrire{#1secno}{\the\secno}\else\relax
\ecrire{#1secno}{\the\secno.\the\parno}\fi\else\relax
\ifnum\parno=0\relax
\ecrire{#1secno}{\the\secno.\the\enonceno}\else\relax
\ecrire{#1secno}{\the\secno.\the\parno.\the\enonceno}\fi\fi} 

\def\ref#1{\ifundefined{#1chapno}\maquereau={??}\else\ifnum
\csname #1chapno\endcsname=\chapno\maquereau={\csname #1secno\endcsname }\else\ifnum
\csname #1chapno\endcsname>64\relax
\maquereau={{\char\csname #1chapno\endcsname }.\csname #1secno\endcsname }\else
\maquereau={\uppercase
\expandafter{\romannumeral\csname #1chapno\endcsname}.\csname #1secno\endcsname }\fi
\fi\fi\the\maquereau} 

%
%


\def\addchaptocontents#1{\relax}
\def\addsectocontents#1{\relax}
\def\addpartocontents#1{\relax}
\def\addbibtocontents#1{\relax}
\def\addintrotocontents#1{\relax}
\def\addinitapptocontents{\relax}

\def\leaderfill{\leaders\hbox to 1em{\hss. \hss}\hfill}
\def\page #1{\unskip\leaders\hbox to 5mm{\hss.\hss}\hfill
\kern -.7\rightskip\rlap{\hbox to \rightskip{\hss #1}}
\smallbreak}

\def\partentry#1#2{\global\advance\chapno by 1
\secno=0
\ifnum\chapno>64\relax
\def\id{\char\the\chapno}\else
\def\id{\uppercase\expandafter{\romannumeral\the\chapno}}\fi
\hbox to 148mm{\id . #1\leaderfill #2}}

\def\sectionentry#1#2{\global\advance\secno by 1
\parno=0
\ifnum\chapno=0\relax
\def\id{\the\secno}\else
\def\id{\the\secno}\fi
\smallbreak
\vbox{\hskip -1.2cm\hsize=146.5mm{\bf\id . #1\leaderfill #2}}\smallbreak}

\def\paragrentry#1#2{\global\advance\parno by 1
\ifnum\chapno=0\relax
\def\id{\the\secno .\the\parno}\else
\def\id{\the\secno .\the\parno}\fi
\smallbreak
\vbox{\hskip -1.4cm\hsize=146.5mm{\qquad\id . #1\leaderfill #2}}\smallbreak}

\def\bibentry#1#2{
\smallbreak
\vbox{\hskip -1.2cm\hsize=146.5mm{\bf #1\leaderfill #2}}\smallbreak}

\def\introentry#1#2{
\smallbreak
\vbox{\hskip -1.2cm\hsize=146.5mm{\bf #1\leaderfill #2}}\smallbreak}

\def\addintrotocontentsbase#1{\write\contents{\noexpand\introentry{#1}{\the\pageno}}}
\def\addchaptocontentsbase#1{\write\contents{\noexpand\partentry{#1}{\the\pageno}}}
\def\addsectocontentsbase#1{\write\contents{\noexpand\sectionentry{#1}{\the\pageno}}}
\def\addpartocontentsbase#1{\write\contents{\noexpand\paragrentry{#1}{\the\pageno}}}
\def\addbibtocontentsbase#1{\write\contents{\noexpand\bibentry{#1}{\the\pageno}}}

\newwrite\contents

\def\initcontents{\immediate\openout\contents=\jobname.tit\relax
\global\let\addintrotocontents=\addintrotocontentsbase
\global\let\addchaptocontents=\addchaptocontentsbase
\global\let\addsectocontents=\addsectocontentsbase
\global\let\addpartocontents=\addpartocontentsbase
\global\let\addbibtocontents=\addbibtocontentsbase
}

\def\makecontents{\chapno=0
\secno=0
\parno=0
\mk\ilestpreferabledecouperici
\centerline{\fourteenbf Table des matières}\nobreak
\vskip 15mm\input\jobname.tit
\chapno=0
\secno=0
\parno=0}




 \rm
\font\tenpc=cmcsc10

\newfam\pcfam
\textfont\pcfam=\tenpc
\def\pc{\fam\pcfam\tenpc}


\font\tencyr=cmbx10
\font\sevencyr=cmbx7
\font\fivecyr=cmbx6

\newfam\cyrfam
\textfont\cyrfam=\tencyr
\scriptfont\cyrfam=\sevencyr
\scriptscriptfont\cyrfam=\fivecyr


\font\tengoth=cmbx10
\font\eightgoth=cmbx8
\font\sevengoth=cmbx7

\font\fivegoth=cmbx5

\newfam\gothfam
\textfont\gothfam=\tengoth
\scriptfont\gothfam=\sevengoth
\scriptscriptfont\gothfam=\fivegoth
\def\goth{\fam\gothfam}


\font\tenbboard=msbm10
 \font\eightbboard=msbm8
 \font\sevenbboard=msbm7

 \newfam\bboardfam
 \textfont\bboardfam=\tenbboard
 \scriptfont\bboardfam=\sevenbboard
 \scriptscriptfont\bboardfam=\sevenbboard
 \def\bbfont{\fam\bboardfam\tenbboard}
 \def\bb#1{{\bbfont #1}}

 
 \font\twelvebf=cmbx10 at 12pt
\font\fourteenbf=cmbx10 at 14pt

\font\twelverm=cmr10 at 12pt 

\font\eightpc=cmcsc8
\font\eightrm=cmr8 \font\eightsl=cmsl8 \font\eightbf=cmbx8
\font\eightsy=cmsy8 \font\eighti=cmmi8 \font\eightit=cmti8
\font\sixrm=cmr6   
\font\sixsy=cmsy6 \font\sixi=cmmi6 
\font\sixi=cmmi6 

\skewchar\eighti='177 \skewchar\eightsy='60
\skewchar\sixi='177 \skewchar\sixsy='60

\def\eightpoint {%
\textfont0=\eightrm \scriptfont0=\sixrm
\scriptscriptfont0=\fiverm \def\rm{\fam0\eightrm}%
\textfont1=\eighti \scriptfont1=\sixi
\scriptscriptfont1=\fivei \def\oldstyle{\fam1\eighti}%
\textfont2=\eightsy \scriptfont2=\sixsy
\scriptscriptfont2=\fivesy 
\textfont\itfam=\eightit \def\it{\fam\itfam\eightit}%
\textfont\slfam=\eightsl \def\sl{\fam\slfam\eightsl}%
\textfont\pcfam=\eightpc \def\pc{\fam\pcfam\eightpc}%
\textfont\bffam=\eightbf 
 \def\bf{\fam\bffam\eightbf}%
\textfont\gothfam=\eightgoth 
\def\goth{\fam\gothfam\eightgoth}%
\textfont\bboardfam=\eightbboard 
\def\bb{\fam\bboardfam\eightbboard}%
\abovedisplayskip=9pt plus 2pt minus 6pt
\belowdisplayskip=\abovedisplayskip
\abovedisplayshortskip=0pt plus 2pt
\belowdisplayshortskip=5pt plus 2pt minus 3pt
\smallskipamount=2pt plus 1pt minus 1pt
\medskipamount=4pt plus 2pt minus 2pt
\bigskipamount=9pt plus 4pt minus 4pt
\setbox\strutbox=\hbox{\vrule height 7pt depth 2pt width 0pt}%
\normalbaselineskip=9pt \normalbaselines
\rm}

\def\tenpoint {%
\textfont0=\tenrm \scriptfont0=\sevenrm
\scriptscriptfont0=\fiverm \def\rm{\fam0\tenrm}%
\textfont1=\teni \scriptfont1=\seveni
\scriptscriptfont1=\fivei \def\oldstyle{\fam1\tenti}%
\textfont2=\tensy \scriptfont2=\sevensy
\scriptscriptfont2=\fivesy 
\textfont\itfam=\tenit \def\it{\fam\itfam\tenit}%
\textfont\slfam=\tensl \def\sl{\fam\slfam\tensl}%
\textfont\bffam=\tenbf 
\def\bf{\fam\bffam\tenbf}%
\textfont\gothfam=\tengoth 
\def\goth{\fam\gothfam\tengoth}%
\textfont\bboardfam=\tenbboard 
\def\bb{\fam\bboardfam\tenbboard}%
\abovedisplayskip=12pt plus 3pt minus 9pt
\belowdisplayskip=\abovedisplayskip
\abovedisplayshortskip=0pt plus 3pt
\belowdisplayshortskip=7pt plus 3pt minus 4pt
\smallskipamount=3pt plus 1pt minus 1pt
\medskipamount=6pt plus 2pt minus 2pt
\bigskipamount=12pt plus 4pt minus 4pt
\setbox\strutbox=\hbox{\vrule height 8.5pt depth 3.5pt width 0pt}%
\normalbaselineskip=12pt \normalbaselines
\rm}


\def\cqfd{\unskip\hfill\hbox{\vrule\vbox to 6pt{\hrule width 4pt
\vfill\hrule}\vrule}}
 

\def\diag#1{\def\normalbaselines{\baselineskip=0pt
                \lineskip=10pt\lineskiplimit=1pt} \matrix{#1}}

\def\hfl#1#2{\smash{\mathop{\hbox to 5mm{\rightarrowfill}}
\limits^{\scriptstyle#1}_{\scriptstyle#2}}}

\def\Hfl#1#2{\smash{\mathop{\hbox to 25mm{\rightarrowfill}}
\limits^{\scriptstyle#1}_{\scriptstyle#2}}}

\def\sec#1{\smash{\mathop{\mathrel<\joinrel\mathrel
{\hbox to 5mm{\dotfill}}}\limits_{\scriptstyle#1}}}
\def\sect#1#2{\smash{\mathop{\hbox to 5mm{\rightarrowfill}}
\limits^{\scriptstyle#1}_{\sec{#2}}}}

\def\ghfl#1#2{\smash{\mathop{\hbox to 5mm{\leftarrowfill}}
\limits^{\scriptstyle#1}_{\scriptstyle#2}}}

\def\rto{\rightarrow}

\def\rpair#1#2{\raise -2pt\vbox{\baselineskip=4pt
\hbox{$\displaystyle\hfl{#1}{}$}
\hbox{$\displaystyle\hfl{}{#2}$}}}
\def\lrpair#1#2{\raise -2pt\vbox{\baselineskip=4pt
\hbox{$\displaystyle\ghfl{#1}{}$}
\hbox{$\displaystyle\hfl{}{#2}$}}}


\def\build#1_#2^#3{\mathrel{\mathop{\kern 0pt#1}\limits_{#2}^{#3}}}


\def\op{ope\-rad}
\def\alg{alge\-bra}

\def\gr{gra\-ded}
\def\vs{vector space}
\def\pl{pre-Lie}
\def\pla{\pl\ \alg}
\def\rtt{roo\-ted tree}

\def\mor{mor\-phism}
\def\iso{iso\-\mor}


\def\ot{\otimes}
\def\oco{\ot\cdots\ot}
\def\nd{\noindent}

\def\thr{theorem}

\def\defi{de\-fi\-ni\-tion}

\def\st{such that}

\def\newfonction#1{\expandafter\def\csname#1\endcsname{\mathop{\rm #1}\nolimits}}
\newfonction{In}
\newfonction{Id}
\newfonction{Tor}
\newfonction{Perm}
\newfonction{Ker}
\newfonction{Coder}
\newfonction{Im}
\newfonction{sgn}
\newfonction{pr}
\newfonction{Hom}
\newfonction{dim}
\newfonction{Pe}
\newfonction{CPL}
\newfonction{Sh}

\def\S{{\rm S}}

\def\Z{{\bb Z}}
\def\SS{{\bb S}}

\def\call#1{\expandafter\def\csname#1\endcsname{{\cal #1}}}

\call{P}
\call{C}
\call{D} 
\call{A}
\call{L}
\call{M}
\call{N}
\call{F}
\call{G}
\call{U}

\def\Lie{\L{\rm ie}}

\def\PL#1{\mathop{\cal{PL}}#1}
\def\RT#1{\mathop{\cal{RT}}#1}
\def\HPL{\mathop{\rm HPL}}




\centerline{\twelvebf PRE-LIE ALGEBRAS AND THE ROOTED TREES OPERAD}
\vskip 0.4cm
\hbox{\hfill
\vbox to 39mm{\vfill
\hbox{\twelverm  Fr\'ed\'eric CHAPOTON\qquad and}
\hbox{\eightpoint Institut de Math\'ematiques,}
\hbox{\eightpoint Equipe d'Analyse alg\'ebrique}
\hbox{\eightpoint 4, place Jussieu}
\hbox{\eightpoint 75252 Paris Cedex 05}
\hbox{\eightpoint chapoton@math.jussieu.fr}
\vfill}
\hskip 1cm
\vbox to 39mm{\vfill
\hbox{\twelverm  Muriel LIVERNET}
\hbox{\eightpoint LAGA, Institut Galil\'ee}
\hbox{\eightpoint Universit\'e Paris 13}
\hbox{\eightpoint Avenue Jean-Baptiste Cl\'ement}
\hbox{\eightpoint 93430 Villetaneuse}
\hbox{\eightpoint livernet@math.univ-paris13.fr}
\vfill}
\hfill}

\vskip -0.2cm


\begingroup
\parindent=0cm
\eightpoint
{\bf Abstract.} A \pla\ is a \vs\ $L$
endowed with a bilinear product $\cdot : L\times L \rto L$
satisfying the relation $(x\cdot y)\cdot z-x\cdot(y\cdot z)=
(x\cdot z)\cdot y-x\cdot(z\cdot y),\ \forall x,y,z \in L$. We give an
explicit combinatorial description in terms of \rtt s of the operad associated to this type
of algebras and prove that it is a Koszul operad.
\par
\endgroup

\vskip -0.4cm
\intro{Introduction}
\vskip -0.3cm
We study here a type of algebra which deserves more attention than it
has been given. People have been using these algebras, under various
names, for a long time. They appeared under the name of left-symmetric
algebras in the work of Vinberg on convex homogeneous cones [V], and so
were dubbed Vinberg algebras in some papers. They also appeared in the
study of affine manifolds, under the name of right-symmetric algebras
[Mat]. We propose to adopt the name of pre-Lie algebras, which
has been used by Gerstenhaber [G]: the Lie bracket involved in the
Gerstenhaber structure on the Hochschild cohomology comes from
a \pla\ structure on the cochains. Besides, rooted trees have shown their
interest in the study of vector fields, numerical analysis
(see e.g. the paper of C. Brouder [B] and the references therein) 
and more recently in quantum field theory [Connes, 
Kreimer]. We define in this paper the underlying \op\ of \pla s in terms
of \rtt s which should shed light on the relationships between these different
topics. The first author is indebted to M. Kontsevich
for a talk about Hochschild complex and Deligne's conjecture which inspired
the link between \rtt s and \pla s.

\medskip
The description of the \op\ defining \pla s  in terms of \rtt s is the subject
of the first section. The \op\ arising here should not be confused
with the structure on \rtt s appearing in [B-V]. The second section is devoted to the definition of
the \op ic  homology of \pla s. Finally, we prove in the third section
that the \op\ associated to \pla s is a Koszul \op. To that end, we prove in fact that a free \pla\ $L$ is a free module over the enveloping \alg\ of the
Lie \alg\ underlying  $L$. Combined with the first section, it gives a new
interpretation of the Hopf \alg\ appearing in the works of 
A. Connes and D. Kreimer [C-K].
\vfill\eject


\section{A description of the \op\ defining \pla s}

This section is devoted to the description of the operad defining
\pla s. We prove in the theorem \ref{thoprtt} that this operad is the operad of \rtt s.
\medskip
We recall briefly some facts about \op s (see [Gi-K],[Ge-J]). An {\it operad} 
$\P$ 
is a sequence of \vs s $\P(n)$, for $n\geq 1$, \st\ $\P(n)$ is a module 
over the
symmetric group $\S_n$, together with composition maps $\gamma:\P(n)\ot\P(i_1)
\oco\P(i_n)\rto\P(i_1+\cdots+i_n)$ satisfying some relations of  
associativity, unitarity and equivariance with respect to the symmetric group,
 called May axioms [May]. Note that giving a 
$\S_n$-module for
all $n$ is equivalent to giving a {\it ${\bb S}$-module}, i.e. a functor from 
the category of (finite sets, bijections) to the category of \vs s. 
Hence an \op\ can be defined by a ${\bb S}$-module $\P$ together 
with composition maps 
$\gamma:\P(I)\ot\P(J_1)\oco\P(J_n)\rto\P(J_1\sqcup\cdots\sqcup J_n)$
where $n$ is the cardinal of $I$. An {\it algebra} over an \op\ $\P$ is
a \vs\ $A$ together with maps $\P(n)\ot A^{\ot n}\rto A$ satisfying 
some relations of associativity, unitarity  and equivariance.


\tit{\pla s}A {\it \pla\ } is a vector space $L$ together with
a bilinear map $\cdot : L\times L\rto L$ satisfying the relation
$$(x\cdot y)\cdot z-x\cdot(y\cdot z)=
(x\cdot z)\cdot y-x\cdot(z\cdot y),\ \forall x,y,z \in L.$$
When the \vs\ $L$ is \gr\ we define a {\it \gr\ \pla\ } by the relation
$$(x\cdot y)\cdot z-x\cdot(y\cdot z)=(-1)^{|y||z|}(
(x\cdot z)\cdot y-x\cdot(z\cdot y)),\ \forall x,y,z \in L.$$
\label{defPLalg}


\PROP Let $(L,\cdot)$ be a pre-Lie algebra. The bracket defined
by 
$$[a,b]=a\cdot b-b\cdot a,\ \forall a,b\in L$$
 endows $L$ with a 
structure of Lie \alg. \fin 
In the sequel $L_{\Lie}$ will denote the Lie \alg\ $(L,[-,-])$.
\label{Lieinduced}

\tit{Examples}M. Gerstenhaber [G] introduced a structure of 
\pla\ on the Hochschild complex of an associative \alg\ $A$ as follows:
denote by $C^m(A,A)$ the space $\Hom(A^{\ot m},A)$ in degree $m-1$;
let $f\in C^m$ and $g\in C^n$, then the product
$$\displaylines{
(f\circ g)(a_1\oco a_{m+n-1})=\cr
\sum_{i=1}^{m} (-1)^{(n-1)(i-1)} 
f(a_1\oco a_{i-1}\ot g(a_i\oco a_{i+n-1})\ot a_{i+n}\oco a_{m+n-1})}$$
satisfies the graded relation defining \pla s.
\medskip
The structure of pre-Lie algebra also appears in the study of affine
structures on manifolds [Mat]. An affine structure on a
$n$-manifold is an atlas whose coordinate changes are in the group of
affine motions of ${\bb R}\sp n$. It can also be given by a linear
connection $\nabla$ whose torsion and curvature vanish. The product
$X\circ Y=-\nabla_Y X$ then defines a structure of pre-Lie algebra on
the set of vector fields, such that the associated Lie bracket is the
usual bracket of vector fields.


\tit{The operad $\PL $}From the \defi\ \ref{defPLalg}, it is clear that
a \pla\ is an \alg\ over a binary quadratic \op, denoted by $\PL $. We 
recall briefly how to construct $\PL $ [Gi-K]. Let $\F$ be the free \op\ generated by the regular representation of $\S_2$. A basis of $\F(n)$, as a \vs, is
given by ``parenthesized products'' on $n$ variables indexed by $\{1,\cdots,n\}$. For 
instance, a basis of $\F(2)$ is given by $(x_1x_2)$ and $(x_2x_1)$, and a 
basis of 
$\F(3)$ is given by $((x_1x_2)x_3)$, $(x_1(x_2x_3))$ and all their 
permutations. Let $R$ be the $S_3$-sub-module of $\F(3)$ generated by the relation 
$r=((x_1x_2)x_3)-(x_1(x_2x_3))-((x_1x_3)x_2)+(x_1(x_3x_2))$. 
Then $\PL =\F/(R)$, where $(R)$ denotes the ideal
of $\F$ generated by $R$. The operadic composition on $\PL $ is induced by the
one on $\F$, given by
$$\gamma :  \F(n)\ot\F(i_1)\oco\F(i_n)\rto \F(i_1+\cdots + i_n)$$   
which assigns to $(\mu,\nu_1,\cdots\nu_n)$ the word obtained by substituting
$\nu_i$ for $x_i \in \mu$. Notice that the concatenation $(\rho\rho')$ is
the particular case of the composition $\gamma((x_1x_2),\rho,\rho')$.
\label{descPL}


\tit{The operad of \rtt s $\RT $}Let $n>0$. A {\it \rtt}  of
degree $n$, or $n$-\rtt, is a non-empty connected 
graph without loops whose vertices are
labelled by the set $[n]=\{1,\cdots, n\}$, together with a 
distinguished element
in this set called the root. Edges of this graph are oriented 
towards the root. We denote by $\RT{(n)}$ the free $\Z$-module generated
by $n$-\rtt s. We can endow $\RT =(\RT{(n)})_{n\geq 1}$
with an operad structure, as explained below.

The action of the symmetric group is the natural one, by permutation of 
indices. 
Let $T$ be a $n$-\rtt; denote
by $\In(T,i)$  the set of incoming edges at the vertex $i$ of $T$. Let
$S$ be a $m$-\rtt. In order to define the operadic composition, 
describing the compositions $\circ_i : \RT{(n)}\times\RT{(m)}\rto
\RT{(n+m-1)}$, for $1\leq i \leq n$ (see e.g. [Lo]) is enough. 
We define the composition
of $T$ and $S$ along the vertex $i$ of $T$ by
$$T\circ_i S=\sum_{f : \In(T,i)\rto [m]} T\circ_i^f S,$$
where $T\circ_i^f S$ is the \rtt\ obtained by substituting the tree $S$ for
the vertex $i$ in $T$~: the outgoing edge of $i$, if exists, becomes the
outgoing edge of the root of $S$; incoming edges of $i$ are grafted on
the vertices  of $S$ following  the map $f$. Then, it is easy to check that 
these compositions endow $\RT$ with a structure of operad. Let us
give an example. A \rtt\ is drawn with its root at the bottom.
Let 
$$T=\vcenter
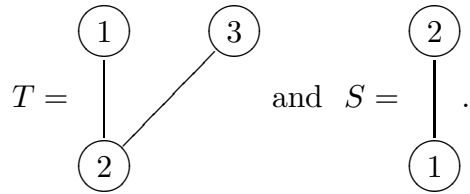
{\xymatrix{*++[o][F-]{1} \ar@{-}[d] & 
*++[o][F-]{3} \ar@{-}[dl]\\
*++[o][F-]{2}}}
\ {\hbox{\rm and}}\ \  S=\vcenter{
\xymatrix{*++[o][F-]{2} \ar@{-}[d] \\
*++[o][F-]{1}}}.$$
By describing maps from $\{1,3\}$ to $\{1,2\}$ and by reindexing the vertices
of $S$ and $T$, one gets
 $$T\circ_2 S=
\!\!\!\!\vcenter{\SelectTips{cm}{}\xymatrix@-1pc{
*++[o][F-]{1}\ar@{-}[d]&*++[o][F-]{4}\ar@{-}[dl]\\
*++[o][F-]{3}\ar@{-}[d]\\
*++[o][F-]{2}\\
{\scriptscriptstyle f: \left\{\matrix{ 1\rto 2 \cr 3\rto 2}\right.}
}}\!\!\!\!\!\!\!\!+\!\! 
\vcenter{\SelectTips{cm}{}\xymatrix@-1pc{
*++[o][F-]{1}\ar@{-}[d]\\
*++[o][F-]{3}\ar@{-}[d]&*++[o][F-]{4}\ar@{-}[dl]\\
*++[o][F-]{2}\\
{\scriptscriptstyle f: \left\{\matrix{ 1\rto 2 \cr 3\rto 1}\right.}
}}\; +\!\!\!\!\!\!
\vcenter{\SelectTips{cm}{}\xymatrix@-1pc{
*++[o][F-]{4}\ar@{-}[d]\\
*++[o][F-]{3}\ar@{-}[d]&*++[o][F-]{1}\ar@{-}[dl]\\
*++[o][F-]{2}\\
{\scriptscriptstyle f: \left\{\matrix{ 1\rto 1 \cr 3\rto 2}\right.} }}\; +\!\!
\raise 1.65cm \vbox{\SelectTips{cm}{}\xymatrix@-1pc{
\\
*++[o][F-]{1}\ar@{-}[dr]&*++[o][F-]{3}\ar@{-}[d]&*++[o][F-]{4}\ar@{-}[dl]\\
&*++[o][F-]{2}\\
&{\scriptscriptstyle f: \left\{\matrix{ 1\rto 1 \cr 3\rto 1}\right.}
}}.$$


\tit{The Poincare series associated to $\RT $}The Poincare series is defined by

$$g_{\RT{}}(x)=\sum_{n\geq 1} \dim(\RT{(n)})(-x)^n/n!$$
Using classical results in combinatorics [W], one obtains that 
$\dim(\RT{(n)})=n^{n-1}$ and that $g_{\RT{}}$ is the inverse map of 
$x\mapsto -xe^{-x}$.


\tit{A product in the operad $\RT{}$}Recall that the operadic composition
for $\RT{}$, $\gamma : \RT{(n)}\ot\RT{(i_1)}\oco
\RT{(i_n)}\rto\RT{(i_1+\cdots + i_n)}$ is given in terms of 
the $\circ_i$ compositions by 
$\gamma(\mu,\nu_1,\cdots,\nu_n)=
(\cdots(\mu\circ_{n}\nu_n)\circ_{n-1}\nu_{n-1})\cdots\circ_1\nu_1).$ 

We may then define a particular composition, which corresponds to the concatenation in
$\PL $-case (see \ref{descPL}): for any \rtt s $T_1$ and $T_2$, let
$$

T_1\star T_2=\gamma(
\vcenter{
\xymatrix@-1.5pc{*++[o][F-]{2} \ar@{-}[d] \\
*++[o][F-]{1}}},T_1,T_2)=\sum_{s\in I_1}
\vcenter{\xymatrix@-1.5pc{
*++[o][F-]{T_2}\ar@{-*{\bullet}}[d]_>>{s}\\
*++[o][F-]{T_1}}}.$$
More explicitly, the operation
$T_1\star T_2$ consists
of grafting the root of $T_2$ on every vertex of $T_1$.
\label{prodRT}


\LEM The following relation holds for any \rtt s $T_j$, $1\leq j\leq 3$:
$$(T_1\star T_2)\star T_3 -T_1\star (T_2\star T_3)= 
(T_1\star T_3)\star T_2- T_1\star (T_3\star T_2).$$\fin
\label{lemmRPL}
\vskip -0.4cm
\dem By computing
$$(T_1\star T_2)\star T_3-T_1\star(T_2\star T_3)=

\sum_{s\in I_1}\sum_{t\in I_2}
\vcenter{\xymatrix@-1.5pc{
*++[o][F-]{T_3}\ar@{-*{\bullet}}[d]_>>{t}\\
*++[o][F-]{T_2}\ar@{-*{\bullet}}[d]_>>{s}\\
*++[o][F-]{T_1}}}+\sum_{s\in I_1}\sum_{s\in I_2}
\vcenter{\xymatrix@-1.5pc{
*++[o][F-]{T_3}\ar@{-*{\bullet}}[d]_>>{t}
&*++[o][F-]{T_2}\ar@{-*{\bullet}}[dl]^>{s}\\
*++[o][F-]{T_1}}}
-\sum_{s\in I_1}\sum_{t\in I_2}
\vcenter{\xymatrix@-1.5pc{
*++[o][F-]{T_3}\ar@{-*{\bullet}}[d]_>>{t}\\
*++[o][F-]{T_2}\ar@{-*{\bullet}}[d]_>>{s}\\
*++[o][F-]{T_1}}}$$
and inverting the roles of $T_2$ and $T_3$, the required equality is 
obtained. \cqfd


\TH The operad $\PL $ defining \pla s is isomorphic to the
\op\ of \rtt s $\RT{}$. \fin
\label{thoprtt}
\dem Firstly, we define an \op ic \mor\ $\Phi : \PL \rto\RT $. 
Since $\PL =\F/(R)$ 
(\ref{descPL}) it is sufficient to define $\Phi$ on $\PL{(2)}=\F(2)$,
then to extend it on $\F$ by the universal property of the free \op\ and
to check that $\Phi(r)=0$. Set 
$$

\Phi((x_1x_2))=\vcenter{\xymatrix@-1.5pc{
*++[o][F-]{2}\ar@{-}[d]\\
*++[o][F-]{1}}}\quad\hbox{\rm and}\quad
\Phi((x_2x_1))=\vcenter{\xymatrix@-1.5pc{
*++[o][F-]{1}\ar@{-}[d]\\
*++[o][F-]{2}}}.$$ 
Hence
$$

\Phi(r)=\left(
\vcenter{\xymatrix@-1.2pc{
*++[o][F-]{3}\ar@{-}[d]\\
*++[o][F-]{2}\ar@{-}[d]\\
*++[o][F-]{1}}}+ 
\vcenter{\xymatrix@-1.2pc{
*++[o][F-]{3}\ar@{-}[d]&*++[o][F-]{2}\ar@{-}[dl]\\
*++[o][F-]{1}}}\right) -
\vcenter{\xymatrix@-1.2pc{
*++[o][F-]{3}\ar@{-}[d]\\
*++[o][F-]{2}\ar@{-}[d]\\
*++[o][F-]{1}}}-
\left(\vcenter{
\xymatrix@-1.2pc{
*++[o][F-]{2}\ar@{-}[d]\\
*++[o][F-]{3}\ar@{-}[d]\\
*++[o][F-]{1}}}+ 
\vcenter{\xymatrix@-1.2pc{
*++[o][F-]{2}\ar@{-}[d]&*++[o][F-]{3}\ar@{-}[dl]\\
*++[o][F-]{1}}}\right)+
\vcenter{\SelectTips{cm}{}\xymatrix@-1pc{
*++[o][F-]{2}\ar@{-}[d]\\
*++[o][F-]{3}\ar@{-}[d]\\
*++[o][F-]{1}}}=0.$$

Remark that, since $\Phi$ is an \op ic \mor, it sends the concatenation
product (see \ref{descPL}) to the product $\star$ defined in \ref{prodRT}.
\tibla
The proof relies on the existence of an inverse $\Psi$ of $\Phi$.
As we explained in the introduction of the section, it is more convenient
for the proof to deal with $\SS$-modules. For a finite set
$I$, a {\it $I$-\rtt} is a \rtt\ whose vertices
are labelled by $I$; its degree is the cardinal of $I$. 
We denote by $\Phi_I : \PL{(I)}\rto\RT{(I)}$ the natural extension of $\Phi$.
\tibla
\nd{\sl Claim. }\debut For any finite set $I$, there exists
a map $\Psi_I : \RT{(I)}\rto\PL{(I)}$ \st\ $\Psi_I\Phi_I=\Id$
and $\Phi_I\Psi_I=\Id$.\fin
\tibla
We prove the claim by induction on $\#I$. If $I=\{i\}$, then it is trivial. 
Assume that the claim is true for any $I$ \st\ $\#I\leq n$ and let $I$ be
a finite set of cardinal $n+1$. Let $T$ be a $I$-rooted tree and $i$ be
its root. Up to a permutation, we can write uniquely
$$T=B(i,T_1,\cdots, T_p)=

\vcenter{
\SelectTips{cm}{}
\xymatrix@-1.5pc{
*++[o][F-]{T_1} \ar@{-}[dr] & *++[o][F-]{T_2}\ar@{-}[d]&
\cdots\cdots&
*++[o][F-]{T_p} \ar@{-}[dll]\\
& *++[o][F-]{i}}},$$
where $T_i$, for $1\leq i\leq p$, is a \rtt\ of degree 
strictly less than $n+1$. Let us define the map $\Psi_I$ by induction on $p$.

If $p=1$, then $T=B(i,T_1)=

\xymatrix@-1.5pc{*++[o][F-]{i}}\star T_1$ and $\Psi_I(T)=(i\Psi_I(T_1))$
is well defined. Moreover $\Phi_I\Psi_I(T)=T$, because $\Phi_I$
sends concatenation to the product $\star$ and because of the induction
hypothesis. For $p\geq 2$, one has 
$$T=B(i,T_2,\cdots T_p)\star T_1-
\sum_{j=2}^{p} B(i,T_2,\cdots,T_j\star T_1,\cdots T_p).$$
Consequently, we may define by induction: 
$$\Psi_I(T)=(\Psi_I(B(i,T_2,\cdots T_p))\Psi_I(T_1))-
\sum_{j=2}^{p} \Psi_I(B(i,T_2,\cdots,T_j\star T_1,\cdots T_p)).$$ 
Moreover, as  in the case $p=1$, $\Phi_I\Psi_I(T)=T$.
\tibla
A priori, since $T$ is uniquely determined only up to a permutation,
this definition depends on the choice of
the edge we ungraft in the tree $T$. Let us prove by induction on $p$
that it is not the case. For $p=0,1$ there is no choice; for $p>1$,
we prove that ungrafting the edge where $T_1$ lies, then
ungrafting edges where $T_2$ lies in the trees involved in the sum,
 gives the same definition of $\Psi_I(T)$
as doing it inverting  $T_1$ and  $T_2$. By ungrafting $T_2$ in the previous
relation, we get
$$\displaylines{
T=(B(i,T_3,\cdots T_p)\star T_2)\star T_1
-\sum_{k=3}^{p}B(i,T_3,\cdots,T_k\star T_2,\cdots, T_p)\star T_1\cr 
-B(i,T_3,\cdots T_p)\star (T_2\star T_1)
+\sum_{j=3}^{p} B(i,T_3,\cdots,T_j\star (T_2\star T_1),\cdots T_p)\cr
- \sum_{j=3}^{p}B(i,T_3,\cdots,T_j\star T_1,\cdots T_p)\star T_2
+\sum_{j=3}^{p}B(i,T_3,\cdots,(T_j\star T_1)\star T_2,\cdots T_p)\cr
+\sum_{j=3}^{p}\sum_{k=3\atop k\not=j}^{p}
B(i,T_3,\cdots,T_k\star T_2,\cdots,T_j\star T_1,\cdots, T_p).}$$
Let $A_{12}=((\Psi_IB(i,T_3,\cdots T_p)\Psi_I(T_2))\Psi_I(T_1))
-(\Psi_IB(i,T_3,\cdots T_p)(\Psi_I (T_2)\Psi_I(T_1)))$ 
and let $A_{21}$ be
the same term with $T_1$ and $T_2$ inverted; since
$A_{21}-A_{12} \in (R)$,  these terms coincide
in $\PL$. It is clear that the terms
$B_{12}^{k}=(\Psi_IB(i,T_3,\cdots,T_k\star T_2,\cdots, T_p)\Psi_I(T_1))
+(\Psi_IB(i,T_3,\cdots,T_k\star T_1,\cdots, T_p)\Psi_I(T_2))$
and $B_{21}^{k}$ coincide, as well as the terms \hfill\break$C_{12}^{jk}=
\Psi_I(B(i,T_3,\cdots,T_k\star T_2,\cdots,T_j\star T_1,\cdots, T_p))$
and $C_{21}^{kj}$; finally, the terms
$D_{12}^j=\Psi_I(B(i,T_3,\cdots,(T_j\star T_1)\star T_2,\cdots T_p))
+\Psi_I(B(i,T_3,\cdots,T_j\star (T_2\star T_1),\cdots T_p)$ and
$D_{21}^{j}$ coincide thanks to lemma \ref{lemmRPL}. 
\tibla
Furthermore $\Psi\Phi=\Id$. In fact, one has $\Psi(T\star T')=
(\Psi(T)\Psi(T'))$. Indeed when $T=
\xymatrix@-2pc{*++[o][F-]{i}}$ the result
comes from the case $p=1$; if $T=B(i,T_1,\cdots, T_p)$,
then $T\star T'=B(i,T',T_1,\cdots,T_p)+\sum_{k=1}^{p}B(i,T_1,\cdots,
T_k\star T',\cdots, T_p)$, and by choosing to ungraft $T'$ in order to 
 define $\Psi$, we get the result. As a consequence, 
let $\mu$ be a word
in $\PL$, hence $\mu$ can be uniquely decomposed in a concatenation
$\mu=(\rho\rho')$; since $\Phi(\mu)=\Phi(\rho)\star\Phi(\rho')$,
then $\Psi\Phi(\mu)=(\Psi\Phi(\rho)\Psi\Phi(\rho'))$ and we can conclude
by induction on the degree of $\mu$. This ends the proof of the
theorem. \cqfd
\tibla
The description of free \pla s is a direct consequence of the previous results.
\tibla
\COR Let $V$ be a vector space. The free \pla\ generated by $V$ is
the \vs\ generated by the \rtt s labelled by a basis of $V$, with the product
$\star$ defined in \ref{prodRT} : let $T_1$ and $T_2$ be two trees labelled
by a basis of $V$ then $T_1\star T_2$ is the sum over the vertices $v$ of 
$T_1$ of trees obtained by linking with an edge the root of $T_2$ to the vertex
$v$ of $T_1$. \fin
\label{freePL}
\smallskip
\noindent For instance, let 
$T_1=\vcenter
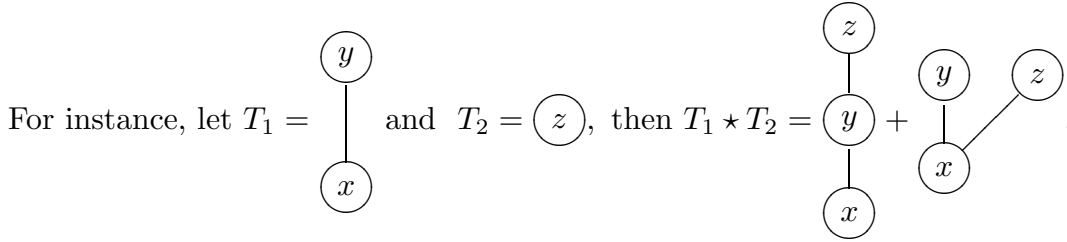
{\xymatrix{*++[o][F-]{y} \ar@{-}[d] \\ 
*++[o][F-]{x}}}
\ {\hbox{\rm and}}\ \  T_2=\vcenter{
\xymatrix{*++[o][F-]{z}}},\ {\hbox{\rm then}}$
$T_1\star T_2=
\vcenter{\SelectTips{cm}{}\xymatrix@-1pc{
*++[o][F-]{z}\ar@{-}[d]\\
*++[o][F-]{y}\ar@{-}[d]\\
*++[o][F-]{x}\\}}+ 
\vcenter{\SelectTips{cm}{}\xymatrix@-1pc{
*++[o][F-]{y}\ar@{-}[d]&*++[o][F-]{z}\ar@{-}[dl]\\
*++[o][F-]{x}\\}}.$

\vskip -0.2cm
\section{Homology of \pla s}

\vskip -0.3cm
Since the homology of \pla s has been  already defined by A. Nijenhuis [N], 
and extended by A. Dzhumadil'daev [D], the aim of this section is to understand
how the operad theory can lead naturally to the definition of the operadic
homology of a type of \alg, which coincides hopefully with some definitions
given earlier.

Following Ginzburg and Kapranov [Gi-K], in order to define the operadic homology
of a \pla, it is necessary to introduce the complex built on the free coalgebra
on the dual \op\ of $\PL$ whose differential is the only coderivation 
induced by the  pre-Lie product.
\medskip
\nd {\bf Notation.} In the sequel, the ground field $K$ will be
of characteristic zero. 

A {\it $(k_1,\cdots, k_p)$-shuffle} is a permutation
$\sigma \in S_{k_1+\cdots +k_p}$ such that $\sigma(1)<\cdots <\sigma(k_1)$,
$\sigma(k_1+1)<\cdots<\sigma(k_1+k_2)$, etc... 
We denote by $\Sh_{k_1,\cdots,k_p}$ the set of 
all $(k_1,\cdots,k_p)$-shuffles. A permutation $\sigma\in S_n$ acts on 
$V^{\ot n}$ by $\sigma\cdot(v_1\oco v_n)=v_{\sigma(1)}\oco v_{\sigma(n)}$, and
in case $V$ is a \gr\ \vs\ we denote by $\epsilon(\sigma,\bar v)$ the sign
appearing in this action. For instance
$\epsilon((12),\bar v)=(-1)^{|v_1||v_2|}$.

Let $V$ be a \gr\ \vs; we denote by $S(V)$ the {\it \gr\ symmetric \alg} 
generated by $V$, that means $S(V)$ is the quotient of the free associative
algebra $T(V)$ generated by $V$ by the ideal generated by $x\ot y-(-1)^{|x||y|}
y\ot x,\ \forall x,y\in V$. If $V$ is concentrated in degree 1, it becomes the
exterior algebra on $V$ and it is usually denoted by $\Lambda(V)$. The
{\it suspension} of $V$, denoted by $sV$ is defined by $(sV)_n=V_{n-1}$.


\PROP The dual \op\ of the pre-Lie \op\ is the \op\ $\Perm$ defined in [Ch]:
a $\Perm$-algebra is a \vs\ $A$ together with a bilinear product
$\cdot : A\times A\rto A$ satisfying the relations
$$\eqalign{
& (a\cdot b)\cdot c= a\cdot (b\cdot c), \cr
& a\cdot b\cdot c=a\cdot c\cdot b,\ \forall a,b,c \in A.}$$ \fin
\label{dualPL}

\dem Recall the definition of the quadratic dual \op\ 
$\P^!$ of an operad $\P$ in
our framework, with the notation of \ref{descPL}. If $\P=\F/(R)$, where
$\F$ is the free \op\ generated by the regular representation of $\S_2$, 
then there is a scalar
product on $\F(3)$ defined by
$$\eqalign{
<i(jk),i(jk)>=&\sgn\pmatrix{1&2&3\cr i&j&k}\cr
<(ij)k,(ij)k)>=&-\sgn\pmatrix{1&2&3\cr i&j&k}}$$
and  $\P^!=\F/(R^\bot)$ where $R^\bot$ is the annihilator 
of $R$ with respect to
this scalar product [Gi-K].

Let $R$ be the $S_3$-sub-module of $\F(3)$ defined in \ref{descPL} and 
$R'$ be the $S_3$-sub-module of $\F(3)$ generated by the relations 
$s=((x_1x_2)x_3)-(x_1(x_2x_3))$ and $t=((x_1x_2)x_3)-((x_1x_3)x_2)$, 
since $<R,R'>=0$, 
$\dim(R)=3$, and $\dim(R')=9$, we can conclude that $R'=R^\bot$.\cqfd

\DEF A {\it graded $\Perm$-coalgebra} $C$ is a positively graded \vs\ equip\-ped 
with a comultiplication $\Delta : C\rto C\ot C$ of degree $0$, satisfying 
the following identities
$$\eqalign{
&(\Id\ot\Delta)\Delta=(\Delta\ot\Id)\Delta, \cr
&(\Id\ot\Delta)\Delta=(\Id\ot T)(\Id\ot\Delta)\Delta}$$
where $T(a\ot b)=(-1)^{|a||b|}b\ot a$.

A {\it coderivation} of $C$ is a linear map $d:C\rto C$ satisfying $\Delta d=
(\Id\ot d)\Delta+(d\ot\Id)\Delta$. The space of all coderivations of $C$
is denoted by $\Coder(C)$.
\medskip
The proof of the following lemma is left to the reader.

\LEM Let $V$ be a reduced graded \vs, i.e. $V_0=0$, then the free
$\Perm$-coalgebra generated by $V$, denoted by $\Perm^c(V)$,
is the \vs\ $V\ot S(V)$ equipped with
the following comultiplication :
$$\eqalign{
&\Delta(v\ot 1)=0, \forall v\in V\cr
&\Delta(v_0\ot  v_1\cdots v_n)=
\sum_{0\leq k\leq n-1\atop\sigma\in \Sh_{k,1,n-1-k}}\!\!\!\!\!\!\!
\epsilon(\sigma,\bar v)\; v_0\ot v_{\sigma(1)}\cdots v_{\sigma(k)}
\ot v_{\sigma(k+1)}\ot
v_{\sigma(k+2)}\cdots v_{\sigma(n)},\cr
&\forall v_0\in V,\ v_1\cdots v_n \in S^n(V).}$$
The space $\Perm^c(V)$ comes with a natural projection onto $V$, denoted
by $\pi$.
\fin

\PROP The following \iso\ of \vs s holds:
\vskip -0.4cm
$$\diag{
\phi: &\Coder(\Perm^c(V))&\rightarrow& \Hom(\Perm^c(V),V)\cr
& d& \mapsto & \pi\circ d.}$$
\vskip -0.4cm
\nd Moreover, there exists a \gr\ \pla\ structure on $\Hom(\Perm^c(V),V)$
such that $\phi$ is an \iso\ of \gr\ Lie \alg s; the Lie bracket
on $\Coder(\Perm^c(V))$ is defined by 
$[d^1,d^2]=d^1d^2-(-1)^{|d^1||d^2|}d^2d^1$. The pre-Lie product $\circ$
is the following one: let $l, m \in\Hom(\Perm^c(V),V)$; 
we set $l=\sum_{i}l_i$ with  $l_i: V\ot S^i(V)\rto V$, 
then 
\vskip -0.2cm
$$\eqalign{
(m\circ l)_n(v_0\ot v_1\cdots v_n)=&\!\!\!\!
\sum_{0\leq i\leq n\atop\sigma\in \Sh_{i,n-i}}\!\!\!\!\!\!
\epsilon(\sigma,\bar v)  
m_{n-i}(l_i(v_0\ot v_{\sigma(1)}\cdots v_{\sigma(i)})\ot 
v_{\sigma(i+1)}\cdots v_{\sigma(n)})\cr
+(-1)^{|v_0||l|}\sum_{0\leq i\leq n-1\atop\sigma\in \Sh_{1,i,n-i-1}}
&\!\!\!\!\!\!
\epsilon(\sigma,\bar v)m_{n-i}(v_0\ot l_{i}(v_{\sigma(1)}\ot 
v_{\sigma(2)}\cdots v_{\sigma(i+1)})
v_{\sigma(i+2)}\cdots v_{\sigma(n)}),}$$
\vskip -0.2cm
\label{propcoder}
\fin

\dem The proof consists essentially of computation. The first part of the proposition relies on the fact that the  map
$\psi : \Hom(\Perm^c(V),V)\rto \Coder(\Perm^c(V))$  defined by
$$\eqalign{
\psi(l)(v_0\ot v_1\cdots v_n)=&\sum_{i=0}^{n}
\sum_{\sigma\in \Sh_{i,n-i}}\!\!\!\!\!\!\epsilon(\sigma,\bar v)\; 
l_i(v_0\ot v_{\sigma(1)}\cdots v_{\sigma(i)})\ot 
v_{\sigma(i+1)}\cdots v_{\sigma(n)}\cr
+(-1)^{|v_0||l|} \sum_{i=0}^{n-1}\sum_{\sigma\in \Sh_{1,i,n-i-1}}&\!\!\!\!\!\!
\epsilon(\sigma,\bar v)\;
v_0\ot l_{i}(v_{\sigma(1)}\ot 
v_{\sigma(2)}\cdots v_{\sigma(i+1)})
v_{\sigma(i+2)}\cdots v_{\sigma(n)},}$$
is the inverse map of $\phi$;  this definition implies also that
$\phi$ is a Lie \alg\ \mor. The fact that the product $\circ$
is a \pla\ \mor\ is a straightforward calculation. \cqfd
\medskip
This proposition yields naturally the definition of a \pla\
up to homotopy.

\DEF A \gr\ \vs\ $V$ is a {\it \pla\ up to homotopy} or a 
{\it $\PL_\infty$-\alg}  if it is
equipped with a map $l\in \Hom(\Perm^c(sV),sV)$ of degree $-1$
\st\ $l\circ l=0$, or equivalently $[l,l]=0$.
\medskip
Furthermore, this proposition leads to the definition 
of the homology of a \pla. Indeed,
let $(L,\cdot)$ be a \pla\ and let $\mu : (sL\ot sL)\rto (sL)$ be the map
of degree $-1$ defined by $\mu(sx\ot sy)=(-1)^{|sx|}s (x\cdot y)$. Hence, by
the definition \ref{defPLalg} of a \pla\ one gets 
$\mu\circ\mu=0$,  then $[\mu,\mu]=0$; thus the coderivation $d$, induced 
by the \iso, satisfies $d^2=0$. The complex so obtained is the one
defining the operadic homology of
a \pla, in the sense of [Gi-K]. Using Koszul sign rules and the \iso\ $sL\simeq e\ot L$, with $e$
a formal element of degree 1, one gets the definition of the homology
of a \pla.
\tibla
\tit{Homology of \pla s}Let $L$ be a \pla. The {\it pre-Lie ho\-mo\-logy} 
of $L$, denoted by $\HPL(L)$ is the homology of the complex $(\CPL_n(L),d)$,
where $\CPL_n(L)=L\ot \Lambda^{n-1}(L)$ and
$$\eqalign{
d(v_0\ot v_1\wedge\cdots\wedge v_n)=&
\sum_{1\leq j\leq n}(-1)^j v_0\cdot v_j \ot v_1\wedge\cdots\wedge \hat v_j
\wedge\cdots\wedge v_n \cr
+\sum_{1\leq i<j\leq n} (-1)^{i+j-1}& v_0\ot [v_i,v_j]
\wedge\cdots\wedge \hat v_i \wedge\cdots\wedge \hat v_j 
\wedge\cdots\wedge v_n,}$$
where the bracket is the Lie bracket in $L$ (see \ref{Lieinduced}).

This complex coincides with the one defined by A. Nijenhuis [N], 
and for a complete definition of the cohomology of a pre-Lie algebra 
with coefficients in a 
representation we will refer to [D]. 

\tibla
\tit{On the link between a \pla\ and its induced Lie \alg} 

Let $(L,\cdot)$ be a \pla\ and denote by $(L_{\Lie},[-,-])$ its induced
Lie algebra (see \ref{Lieinduced}).
Then the relation 
defining the \pla\ structure implies that
$L$ is a right
module over $L_{\Lie}$ via the action
\vskip -0.4cm
$$\diag{
L\times L_{\Lie}& \rto& L\cr
(v,g)&\mapsto & v\cdot g.}$$
\vskip -0.4cm
\nd Hence $L$ is a right module over the enveloping algebra $\U(L_{\Lie})$
of $L_{\Lie}$,  with the usual definition
$l\cdot(a_1\oco a_n)=(\cdots(l\cdot a_1)\cdot a_2)\cdots \cdot a_n).$

As a consequence, there is a nice interpretation of the pre-Lie 
homology of $L$ in terms of the 
Chevalley-Eilenberg homology of $L_{\Lie}$ with coefficients in $L$;
one has the \iso s:
 $${\HPL}_{n+1}(L)\simeq H_n^{CE}(L_{\Lie},L)
\simeq\Tor_n^{\U(L_{\Lie})}(L,K).$$
\label{thhomPL}
\vskip -0.4cm

\section{Koszulness of the \op\ defining \pla s}

The aim of this section is to prove that the \op\ $\PL $ is a Koszul \op. As
explained in [Gi-K], it is enough to prove that for any free \pla\ $L$, its
ho\-mo\-logy $\HPL(L)$ is concentrated in degree $1$.
In fact, the main point is that a free \pla\ $L$ is a free 
right $\U(L_{\Lie})$-module (\thr\ \ref{thFree}). Then the Koszulness of
the \op\ follows with the help of remark \ref{thhomPL}.

Before proving \thr\ \ref{thFree}, we would like to point out some
interesting remarks on free \pla s.

\tit{A link with the Connes-Kreimer hopf algebra}As a Lie algebra, 
the free pre-Lie algebra  on a single generator (see \ref{freePL}) has
already appeared in the work of Alain Connes and Dirk Kreimer on the
combinatorics of renormalization. They consider a commutative Hopf
algebra of polynomials in rooted trees. By the Milnor-Moore theorem,
the dual Hopf algebra is the universal enveloping algebra of some Lie
algebra, which they calculated in [C-K]. This Lie algebra has
a basis indexed by rooted trees, and one can check that the bracket is
the same as the one induced by the pre-Lie structure of the free
pre-Lie algebra.

\LEM Let $L$ be a \pla, $V$ a \vs\ and $\sigma : V\rto L$ a \mor. The right
 $\U(L_{\Lie})$-module $V\ot \U(L_{\Lie})$ can be equipped with a structure
of \pla\ \st\ the map $\tilde\sigma : V\ot \U(L_{\Lie}) \rto L$ defined
by $\tilde\sigma(v\ot u)=\sigma(v)\star u$, where the action of $\U(L_{\Lie})$ on $L$
is denoted by $\star$, becomes a \mor\ of \pla s.\fin
\label{keylemma}

\dem An element in $V\ot U(L_{\Lie})$ will be denoted by $(v,u)$.
The product of \pla\ on $V\ot \U(L_{\Lie})$ is defined as follows
$$(v, u)* (v', u')=(v, u\ot(\sigma(v')\star u')),\ \forall v,v'\in V, u,u'\in
\U(L_{\Lie}).$$
Let us check the relation $R=(A* B)* C-A*(B* C)-
(A* C)* B+A*(C* B)=0$ in $V\ot \U(L_{\Lie})$.
Let $A=(v, u),\ B=(v', u'),\ C=(v'', u'')$.
Then 
$$(A* B)* C-A*(B* C)=(v, u\ot\sigma(v')\star u'\ot \sigma(v'')\star u'')
-(v, u\ot \sigma(v')\star(u'\ot\sigma(v'')\star u'')).$$
But, since $L$ is a right $\U(L_{\Lie})$-module, 
$\sigma(v')\star(u'\ot\sigma(v'')\star u'')=(\sigma(v')\star u')\star
(\sigma(v'')\star u'').$ 
Let $\alpha=\sigma(v')\star u'$ and 
$\beta=\sigma(v'')\star u''$. Since $\alpha$ and $\beta$ lie in $L$ 
the action
$\alpha\star\beta$ coincides with the pre-Lie product in $L$, thus
$\alpha\star\beta-\beta\star\alpha=[\alpha,\beta] \in L_{\Lie}$.
As a consequence
$$R=(v, u\ot (\alpha\ot\beta-\beta\ot\alpha-[\alpha,\beta]))=0.$$
Then it is clear that $\tilde\sigma$ is a \mor\ of \pla s.\cqfd

\TH Let $L$ be a free \pla\ generated by a \vs\ $V$. Then there
is an \iso\ of right $\U(L_{\Lie})$-module
$$L\simeq V\ot \U(L_{\Lie}).$$ \fin
\label{thFree}
\vskip -0.4cm
\dem Let $\sigma$ be the canonical \mor\ from $V$ to $L$,
let $U= \U(L_{\Lie})$ and let $\tau$ be the \mor\ from $V$ to $V\ot U$ \st\ 
$\tau(v)=v\ot 1$. The
product in $L$ is
denoted by $\star$ as well as the action of $U$ on $L$. 
\tibla
By the universal property of the free right $U$-module $V\ot U$, there is a unique
right $U$-module \mor\ $\psi : V\ot U\rto L$ \st\ $\psi\tau=\sigma$. Indeed
$\psi$ is the \mor
\vskip -0.4cm
$$\diag{
\psi : & V\ot U & \rightarrow & L \cr
&v\ot u &\mapsto &\sigma(v)\star u}$$
\vskip -0.4cm
\nd and by virtue of lemma \ref{keylemma}, it
is a \pla s \mor\ for the product  on $V\ot U$ given by
$(v, u)* (v', u')=(v, u\ot(\sigma(v')\star u'))$. Furthermore,
the universal property of the free \pla\ $L$ implies that 
there is a unique \pla s \mor\ $\phi : L\rto V\ot U$ \st\ $\phi\sigma=\tau$.
As a consequence, the fact that $\psi\phi\sigma=\sigma$ and that 
$\psi\phi :L\rto L$
is a \pla\ \mor, implies $\psi\phi=\Id$.
\tibla
In order to conclude it is sufficient to prove that $\phi$ is a right 
$U$-module
\mor, because it implies that $\phi\psi=\Id$. 
The action of $U$ on $V\ot U$ is the concatenation,
denoted by $(v\ot u)\cdot (a_1\oco a_n)=(v, u\ot a_1\oco a_n)$.
Since $U$ is generated by $L$, it is sufficient
to prove that $\phi(x\star y)=\phi(x)\cdot y,\ \forall x,y \in L.$ 
Now $\phi$ is a \pla s \mor, then
$\phi(x\star y)=\phi(x)*\phi(y)$. We set $\phi(x)=\sum_i (v_i, u_i) 
\in V\ot U$
and $\phi(y)=\sum_j (w_j,r_j) \in V\ot U$. Therefore
$$\displaystyle{
\eqalign{
\phi(x)*\phi(y)&=\sum_{i,j}(v_i, u_i\ot \sigma(w_j)\star r_j)\cr
&=\sum_{i,j}(v_i, u_i\ot \psi(w_j\ot r_j))\cr
&=\sum_i (v_i, u_i\ot \psi\phi(y)).}}$$
But $\psi\phi=\Id$, hence $\phi(x)*\phi(y)=\phi(x)\cdot y.$  \cqfd

\TH The \op\ defining \pla s is a Koszul \op.\fin

\dem Let $L$ be the free \pla\ generated by the \vs\ $V$. 
By virtue of \thr\ \ref{thFree}, $L$
is the free right $U(L_{\Lie})$-module on $V$, hence by virtue of remark
\ref{thhomPL}
$$\qquad\qquad\quad {\HPL}_{n}(L)=\Tor_{n-1}^{U(L_{\Lie})}(V\ot U(L_{\Lie}),K)=
\left\{\matrix{
V & \hbox{\rm if}\  n=1,&\qquad\qquad\qquad \cr
0 & \hbox{\rm if not.}\hfill& \cqfd}\right.$$

\medskip
\nd {\bf Acknowledgements.} 
The authors would like to thank  Patrick Polo for his useful 
comments on the third section, James Stasheff and Martin Markl for their interest and suggestions. 
\tibla


\def\colgauche#1{\leavevmode
\item{\vtop to 0pt{\hsize=\parindent\parindent=0pt\leftskip=.5em #1\vss}}}

\biblio{References}
{\eightpoint\parindent=1.2cm
\colgauche{[B-V]}{\pc J.M. Boardman, R.M. Vogt}, ``Homotopy invariant
algebraic structures on topological spaces'', LNM 347, Springer-Verlag, 1973. 
\colgauche{[B]}{\pc C. Brouder}, {\it Runge-Kutta methods and
renormalization}, preprint hep-th/9904014.
\colgauche{[Ch]}{\pc F. Chapoton}, {\it Un endofoncteur de la cat\'egorie
des op\'erades}, preprint.
\colgauche{[C-K]}{\pc A. Connes, D. Kreimer}, {\it Hopf algebras, 
renormalization and noncommutative geometry}, Comm. Math. Phys. 
{\bf 199} (1998), 203-242.
\colgauche{[D]}{\pc A. Dzhumadil'daev}, 
{\it Cohomologies and deformations of right-symmetric algebras}, 
J. Math. Sci. {\bf 93} (1999), 836-876.
\colgauche{[G]}{\pc M. Gerstenhaber} {\it The cohomology structure of 
an associative ring}, Ann. of Math. {\bf 78} (1963), 267-288. 
\colgauche{[Ge-J]}{\pc E. Getzler, J.D.S. Jones}, {\it Operads, homotopy algebra and iterated 
integrals for double loop spaces}, prepublication, 1994.\hfill\break Internet :
 http://xxx.lanl.gov/abs/hep-th/9403055.
\colgauche{[Gi-K]}{\pc V. Ginzburg, M. Kapranov}, {\it Koszul duality for 
operads}, Duke J.Math. (1){\bf 76} (1994), 203-272.
\colgauche{[Lo]}{\pc J.-L. Loday}, {\it La renaissance des op\'erades}, in 
``S\'eminaire Bourbaki, 1994-1995", Ast\'erisque {\bf 237} (1996), 47-74.
\colgauche{[Mat]}{\pc Y. Matsushima}, {\it Affine structures on 
complex manifolds}, Osaka J. Math. {\bf 5} (1968), 215-222.
\colgauche{[May]} {\pc J.P. May} ``The geometry of iterated 
loop spaces'', LNM 271, Springer-Verlag, 1972. 
\colgauche{[N]} {\pc A. Nijenhuis}, {\it On a class of common properties 
of some different types of algebras, I- II} Nieuw Arch. Wisk. {\bf 17} (1969), 
17-46, 87-108.
\colgauche{[V]} {\pc E.B. Vinberg} {\it The theory of homogeneous convex 
cones} (Russian) Trudy Moskov. Mat. Obshch. {\bf 12} (1963) 303-358.
Transl. Moscow Math. Soc. {\bf 12} (1963) 340-403. 
\colgauche{[W]} {\pc H.S. Wilf} ``Generatingfunctionology'' Second edition.
 Academic Press, Inc., Boston, MA, 1994.
\par}

\end